\documentstyle{amsppt}
\magnification =1200

\topmatter
\title Non Dentable Sets in Banach Spaces With Separable Dual
\endtitle

\author Spiros A. Argyros and Irene Deliyanni\\
(Herakleion Crete)\endauthor

\abstract A non RNP Banach space E is constructed such that $E^{*}$
is separable and RNP is equivalent to PCP on the subsets of E.
\endabstract

\endtopmatter

\document

\vskip .1in

The problem of the equivalence of the Radon-Nikodym Property (RNP) and
the Krein Milman Property (KMP) remains open for Banach spaces as
well as for closed convex sets.  A step forward has been made by
Schachermayer's Theorem [S].  That result states that the two
properties are equivalent on strongly regular sets.  Rosenthal, [R],
has shown that every non-RNP strongly regular closed convex set
contains a non-dentable subset on which the norm and weak topologies
coincide.  In a previous paper ([A-D]) we proved that every non RNP
closed convex contains a subset with a martigale coordination.
Furthermore we established the $P\alpha\ell$-representation 
for several cases.  The
remaining open case in the equivalence of RNP and KMP is that of
B-spaces or closed convex sets where RNP is equivalent to PCP in
their subsets.  Typical example for a such structure are
the subsets of $L^1(0,1).$  H. Rosenthal raised the question
if this could occur when the dual of the space is separable.  W.
James $([J_2])$ also posed a similar problem.  The aim of the present
paper is to give an example of a Banach space E with separable dual
failing RNP, and RNP is equivalent to PCP on its subsets.  As
consequence we get that E does not contain $c_o(\Bbb N)$ isomorphically and
hence it does not embed into a Banach space with an unconditional
skipped F.D.D.  On the other hand E semiembeds into a Banach space
with an unconditional basis.  The last property allows us to conclude
that every closed convex non-RNP subset of E contains a closed
non-dentable set with a $P\alpha\ell$-representation.  We recall that a closed
set K has a $P\alpha\ell$-representation if there is an affine, onto, one to one
continuous map from the atomless probability measures on [0,1] to
the set K.  In particular RNP is equivalent to KMP on the subsets
of E.  The space E is realized by applying the Davis-Figiel-Johnson-
Pelczynski factorization method to a convex symmetric set W of a
Banach space $E_u$ constructed in this paper.  Finally as a consequence
of the methods used in the proofs of the example we obtain that every
separable B-space X such that $X^{**}/X$ is isomorphic to $\ell^1(
\Gamma )$ has RNP.
\vskip .1in
We thank H. Rosenthal and T. Odell for some useful discussions
related to the problem studied in the present paper.  We also thank
the Department of Mathematics of Oklahoma State University for its
technical support.
\vskip .1in
We start with some definitions, notations and results necessary for
our constructions.

A closed convex bounded set K is said to be $\delta$-non dentable, $
\delta >0$, if
every slice of K has diameter greater than $\delta$.  A closed convex set
has RNP if it contains no $\delta$-non dentable set.  A closed K subset of a
B-space has the P.C.P. if for every subset L of K and for all $\varepsilon 
>0$ there
exists a relatively weakly open neibhd of L with diameter less than
$\varepsilon$.  It is well known that RNP implies 
P.C.P, but the converse fails [B-R].

In the sequel $\Cal D$ denotes the dyadic tree namely the set of all finite
sequences of the for $a=\{0,\varepsilon_1,...,\varepsilon_n\}$ with $
\varepsilon_i=0\text{ or }1$.  For $a$ in $\Cal D$
the length of $\alpha$ is denoted by $|a|$.  A natural order is induced on $
\Cal D$,
that is $a\prec\beta$ if the sequence $a$ is an initial segment of the sequence
$\beta$.  Two elements $a,\text{ }\beta$ of $\Cal D$ 
are called \underbar{{\bf incomparable}} if they
are imcomparable in the above defined order.  We notice, for later
use, that each $a$ in $\Cal D$ determines a unique basic clopen subset $
V_a$ in
Cantor's group $\{0,1\}^{\Bbb N}$ and $a,\text{ }\beta$ are imcomparable if $
V_a\cap V_{\beta}=\emptyset$.

A basic ingredient in the definition of the space E is Tsirelson's
norm as it is defined in [F-J].  We recall that the norm of this space
satisfies the following implicit fixed point property.

For $\displaystyle{x=\sum^m_{\kappa =1}\lambda_{\kappa}t_{\kappa}}$
$$||\sum^m_{K=1}\lambda_{\kappa}t_{\kappa}||_T=max\{\max_{\kappa}|
\lambda_{\kappa}|,\frac 12sup\sum^n_{j=1}||E_jx||_T\}$$
where the ``$\sup$'' is taken over all choices
$$m<E_1<E_2<...<E_n$$
$E_1,...,E_n$ is an increasing sequence of intervals in the set of natural
numbers and $E_jx$ is the natural projection of x in the space
generated by vectors of the basis $\{t_k:_k\in E_j$ $\}$  Tsirelson's space is a
reflexive Banach space with an unconditional basis not containing any
$\ell^p$ for $1<p<\infty$.
\vskip .1in
{\bf 1.a  The space $\bold E_{\bold u}$
}\vskip .1in

The space $E_u$ will be defined to have an unconditional basis indexed by
the dyadic tree $\Cal D$ and denoted by $(e_a)_{a\in \Cal D}$.  For a sequence of reals
$(\lambda_{\alpha})_{\alpha\in \Cal D}$ which is eventually zero we define
$$\align&||\sum_{a\in \Cal D}\lambda_{\alpha}e_{\alpha}||=\sup \{|
|\sum^{\ell}_{i=1}\lambda_{a_i}t_{k_i}||_T:\{\alpha_i\}^{\ell}_{i=
1}\text{ are incomparable,}\\
&|a_i|=\kappa_i,\text{ }\kappa_1<\kappa_2<...<\kappa_{\ell}\}.\\
\endalign$$
It is clear that $(e_a)_{\alpha\in \Cal D}$ is an unconditional basis for the space $
E_u$ defined
by the above norm.

Next we verify certain properties of the space $E_u$.

\proclaim{1.1  Proposition}The dual of the space $E_u$ is separable.\endproclaim

\demo{Proof}The spare $E_u$ has an unconditional basis hence it is
enough to show that $\ell^1$ does not embed into $E_u$ $[J_1]$.\enddemo

Suppose, on the contrary, that $\ell^1$ embeds into $E_u$.  Then, by standard
arguments, we can find $\ell_1<\ell_2<...<\ell_k<...$ an increasing sequence of
natural numbers and $\{x_k\}^{\infty}_{\kappa =1}$ a normalized sequence in $
E_u$ equivalent
to the usual basis of $\ell^1$, and
$$x_{\kappa}=\sum_{\ell_{\kappa}<|\alpha |<\ell_{\kappa +1}}\lambda_{
\alpha}e_{\alpha}$$

The definition of the norm of $E_u$ and elementary properties of
Tsirelson's norm show that
$$||\sum^m_{k=1}\mu_kx_k||\leq ||\sum^m_{k=1}\text{$_{\kappa}^{\mu}$}
t_{\ell_{k+1}}||_T$$
so $\{t_{\ell_k}\}^{\infty}_{k=2}$ is equivalent to the basis $\ell^
1$.  This contradicts to the
reflexivity of T.
\qed$ $
\vskip .1in

A consequence of the above Proposition is that the basis $(e_{\alpha}
)_{\alpha\in \Cal D}$ is
shrinking.  Therefore every $x^{**}$ in $E_u^{**}$ has a unique representation
as
$$x^{**}=w^{*}\lim_{n\to\infty}\sum_{|\alpha |\leq n}\lambda_{\alpha}
e_{\alpha}:=w^{*}-\sum_{\alpha\in \Cal D}\lambda_{\alpha}e_{\alpha}$$
and $\lambda_{\alpha}=<x^{**},e_{\alpha}$ $^{*}$ $>$.

We define the {\bf support} of $x^{**}$, denoted by supp $x^{**}$, to be the set
$$\{\alpha\in \Cal D:<x^{**},e_{\alpha}^{*}>\neq 0\}.$$
\vskip .1in

\proclaim{1.2  Lemma}Given $x^{**}_1$ $,...,x_{\kappa}^{**}$ in $E_
u^{**}$ such that there are
$a_1,...,a_k$ incomparable elements of $\Cal D$ so that supp $x^{*
*}_i$ is contained in of
$W_{\alpha_i}=\{\beta\in \Cal D:\beta\prec a_i\text{ or }a\prec\beta 
\}$.\endproclaim\ 
Then
$$d(x_1^{**}+..+x^{**}_{\kappa},E_u)\geq\frac 12\sum^{\kappa}_{i=1}
d(x^{**}_1,E_u)$$
\vskip .1in
\demo{Proof}For $n<m$ we define
$$P_{[n,m]}(x^{**})=\sum_{n\leq |a|\leq m}\lambda_{\alpha}e_{\alpha}$$
and
$$P_{[n,\infty ]}(x^{**})=\sum_{n\leq |\alpha |}\lambda_{\alpha}e_{
\alpha}$$
where  $\lambda_{\alpha}=<x^{**},e^{*}_a$ $>$.\enddemo

Using this notation we have
$$d(x^{**},E_u)=\lim_{n\to\infty}||P_{[n,\infty ]}(x^{**})||$$
and
$$||P_{[n,\infty ]}(x^{**})||=\lim_{m\to\infty}||P_{[n,m]}(x^{**})
||$$
\vskip .1in
To establish the result it is enough to show that for $\epsilon >0$ there
exists $n(\varepsilon )$ such that for all $m>n(\varepsilon )$
$$||P_{[m,\infty ]}(\sum^k_{i=1}x^{**}_i)||\geq\frac 12\sum^k_{i=1}
d(x^{**}_i,E_u)-\epsilon .$$
\vskip .1in
Actually $n(\varepsilon )=max\{\kappa ,|a_1|,...|a_{\kappa}|\}$.

Choose any $m>n(\varepsilon )$.  Inductively we define $\{q_i,\ell_
i]^{\kappa}_{i=1}$ such that
$$m<q_1<\ell_1<..<q_k<\ell_k$$
and $||P_{[q_i,\varrho_i]}(x^{**}_i)||>d(x^{**}_i,E_u)-\frac {\in}{
2^i}.$

For each $1\leq i\leq k$ there is a set $\{\beta^i_j\text{ }1\leq 
j\leq s(i)\}$ of incomparable
elements of $\Cal D$ such that $q_i\leq |\beta^i_j|\leq\ell_i$ and
$$||P_{[q_i,\ell_i]}(x^{**}_i)||=||\sum^{s(i)}_{j=1}\lambda_{\beta^
i_j}t_{|\beta^i_j|}||_{_T}.$$
\vskip .in
Notice that $a_i\prec\beta^i_j$ for all $j=1...s(i)$.

Observe that $\displaystyle{\cup_{1\leq i\leq k}}$ $\{\beta^i_j:1\leq 
j\leq s(i)\}$ consists of pairwise
imcomparable elements.  So
$$||P_{[m,\infty ]}(\sum^{\kappa}_{i=1}x^{**}_i)||\geq ||P_{[m,\ell_
k]}(\sum^{\kappa}_{i=1}x^{**}_i)||\geq$$
$$||\sum^{\kappa}_{i=1}\sum^{s(i)}_{j=1}\lambda_{\beta^i_j}t_{|\beta^
i_j|}||_{_T}\geq\frac 12\sum^{\kappa}_{i=1}||\sum^{s(i)}_{j=1}\lambda_{
\beta^i_j}t_{|\beta^i_j|}||_{_T}$$
$$=\frac 12\sum^{\kappa}_{i=1}||P_{[q_i,\ell_i]}(x^{**}_i)||\geq\frac 
12\sum^{\kappa}_{i=1}d(x^{**}_i,E_u)-\epsilon\text{\qed}$$
\vskip .1in
Consider the following closed convex subset of the unit ball of $E_
u$
$$K=\{x\in E_u:x=\sum^{\infty}_{n=0}\sum_{|\alpha |=n}\lambda_{\alpha}
e_{\alpha},\lambda_0=1,\lambda_a\geq 0,\lambda_{\alpha}=\lambda_{(
\alpha ,0)}+\lambda_{(\alpha ,1)}\}.$$
It is easily verified that K is the closed convex hull of a $\frac 
12$-tree
$(d_{\alpha})_{\alpha\in \Cal D}$ where for every $a$ in $\Cal D$ $
d_{\alpha}$ is defined by the conditions
$e^{*}_{\alpha}$ $(d_{\alpha})=1,e^{*}_{(\beta ,0)}(d_{\alpha})=e^{
*}_{(\beta ,1)}(d_{\alpha})=\frac 12e^{*}_{\beta}(d_{\alpha})\text{ and }
d_{\alpha}\in K$.

We set $W=co(K\cup -K)$ and we denote by $\tilde {W}$ its $w^{*}$ closure in $
E^{**}_u$.
Notice that $x^{**}\in\tilde {W}$ if $|e^{*}_{\alpha}(x^{**})|\leq 
1,e^{*}_{(\alpha ,0)}(x^{**})+e^{*}_{(a,1)}(x^{**})=e^{*}_{\alpha}
(x^{**})$
for all $a$ in $\Cal D$.  Hence we could define a map
$$T:M_1(\{0,1\}^{\Bbb N})\rightarrow\tilde {W}$$
with the rule
$$T(\mu )=w^{*}-\sum_{\alpha\in \Cal D}\mu (V_{\alpha})e_a$$
where $V_{\alpha}=\{\gamma\in \{0,1\}^{\Bbb N}:\gamma\restriction 
|\alpha |=a\}$\newline 
Clearly T is one to one and onto.  Furthermore
$$\text{$||T(\mu )||\leq sup\{\sum^k_{i=1}|\mu (V_{\alpha_i})|:\{\alpha_
i\}^k_{i=1}$}\text{ incomparable $\}=||\mu ||.$}$$
Hence T is extended to a bounded linear operator from $M(\{0,1\}^{
\Bbb N})$ onto
the linear span of $\tilde {W}$ denoted by $<\tilde {W}>$.
\vskip .1in

{\bf 1.b  The Space E
}\vskip .1in
The space E is the result of the application of Davis-Figiel-Johnson-Pelczynski 
[D] factorization method to the set W defined above.

We give the precise definition and certain properties of the space E.
$$E=\{y\in E_u:|||y|||=(\sum^{\infty}_{n=1}||y||_n^2)^{\frac 12}<\infty 
\}$$
Here $||.||_n$ denotes the Minkowski's gauge of the set $2^nW+\frac 
1{2^n}B_{E_u}$.

Let $J:E\rightarrow E_u$ be the natural injection.  We notice that $
J[B_E]$
contains the set W; hence E fails RNP.

The operator $J$ satisfies the following properties.
\vskip .1in

P.1.:  $J^{**}:E^{**}\rightarrow E^{**}_u$ is one to one and $J^{*
*}[E^{**}]\cap E_u=J[E].$

As consequence of this property $E^{*}$ is separable.
\vskip .1in

P.2.:  $J$ is a weak to weak homeomorphism on the bounded subsets of
E.  This is a consequence of P.1 and it implies that $J$[L] is closed
for all L, closed convex bounded subsets of E.  In particular $J$ is a
semiembedding.
\vskip .1in

P.3.:  Let L be a closed convex bounded subset of E failing RNP{\bf .}  Then
$J$[L] is non RNP.  If not, $J$[L] is an RNP set, hence for any L-valued
operator $S:L^1\rightarrow E$ the operator $J$o$S$ 
is representable by a function $
\varphi$ in
$L_{J[L]}^{\infty}$.  Then the function $\Psi =J^{-1}\varphi$ 
represents the operator S and K is
RNP.  $([B-R])$
\vskip .1in

P.4.:  If L is bounded subset of E and $J$[L] fails P.C.P. then L fails
P.C.P.

Indeed, for $\{y_n\}^{\infty}_{n=1}$, $y$ in L such that $\displaystyle{\text{$
J(y_n)_{\to}^wJ(y)$}}$ and
$||J(y_n)-J(y)||>\delta >0$ P.2. ensures that $\displaystyle{\text{$
y_n\text{$_{\to}^w$}y$}}$ and also 
$\displaystyle{\text{$||y_n-y||>\frac {\delta}{||J||}$}}$.  Hence y is not a point of continuity.
\vskip .1in

P.5.:  $J^{**}[E^{**}]\subseteq\overline {<\tilde {W}>}$

For this, notice that $B_{E^{**}}\subset 2^n\tilde {W}+\frac 1{2^n}
B_{E_u^{**}}$ hence
$$J^{**}[B_{E^{**}}]\subseteq\cap_n(2^n\tilde {W}+\frac 1{2^n}B_{E^{
**}_u})\subseteq\overline {<\tilde {W}>}.$$
\vskip .1in

We proceed to the proof of the main property of the space E.

\vskip .1in
\proclaim{1.3 Proposition}Let K be a closed, convex, bounded, non RNP
subset of E.  Then K fails P.C.P.\endproclaim
\vskip .1in
\demo{Proof}Property 3, mentioned before, ensures that $J$[K] is non
RNP closed subset of $E_u$.  Hence for some $\delta >0$ there exists a convex
closed L subset of $J$[K] which is $\delta$-nondentable.  Our goal is to show
that every weak neighd in L has diameter greter than $\frac {\delta}{
256}$.  By a
result due to Bourgain [B] it is enough to show that for every
$S_1,S_2,...,S_n$ slices of $\tilde {L}$ there exists $x^{**}_i$ in $
S_i\text{ }i=1,2,...,n$ such that
for all $\displaystyle{(\lambda_i)_{i=1}^n\in \Bbb R^n_{+}}$ 
\hskip0.20in\hskip0.20in$\displaystyle{
\sum^n_{i=1}\lambda_i=1}$
$$d(\sum^n_{i=1}\lambda_ix^{**}_i,E_u)>\frac {\delta}{256}$$
\enddemo

Given $S_1,S_2,...,S_n$ slices of $\tilde {L}$.  Using Lemma 2.7 from [R] we choose
$(x_{\xi ,i}^{**}$ $)_{\xi >\omega_1}$ an uncountable subset of $S_
i$ such that
$$d(x_{\xi ,i}^{**}-x_{\zeta_ii}^{**},E_u)>\frac {3\delta}8\text{ for }
\xi\neq\zeta .$$
\vskip .1in
Recall that $\tilde {L}$ is a subset of $J^{**}[E^{**}]\subset\overline {
<\tilde {W}>}$ and that $T[M\{0,1\}^{\Bbb N}]$ is
norm dense into $\overline {<\tilde {W}>}$.  Hence there are $(\mu_{
\xi ,i})_{\xi <\omega_1,i\leq n}$ such that
$$||T\mu_{\xi ,i}-x_{\xi ,i}^{**}||<\frac {\delta}{256}$$
\vskip .1in

Also, it is known that $\displaystyle{M(\{0,1\}^{\Bbb N})=(\sum_{\gamma 
<2^{\omega}}\oplus L^1(\lambda\gamma ))_1}$ \newline 
where $\{\lambda_{\gamma}\}_{\gamma <2^{\omega}}$ are 
pairwise singular probability measures on $
\{0,1\},^{\Bbb N}$ 
and \linebreak$L^1(\lambda_{\gamma})=L^1[0,1]\text{ or }L^1(\lambda_{
\gamma})=\Bbb R$.

Therefore
$$\mu_{\xi ,i}=\sum_{\gamma <2^{\omega}}\frac {d\mu_{\xi ,i}}{d\lambda_{
\gamma}}$$
where the sum is taken in $\ell_1$-norm.

Choose $F_{\xi ,i}$ finite subset of $2^{\omega}$ so that the measure
$\displaystyle{\mu^{\prime}_{\xi ,i}=\sum_{\gamma\in F_{\xi ,i}}\frac {
d\mu_{\xi ,i}}{d\lambda_{\gamma}}}$ satisfies
$$||T\mu^{\prime}_{\xi ,i}-x_{\xi ,i}^{**}||<\frac {\delta}{256}\qquad
\qquad\qquad\qquad\qquad (1)$$
In particular for $\xi\neq\zeta$ we get
$$d(T\mu^{\prime}_{\xi ,i}-T\mu^{\prime}_{\zeta ,i},E_u)>\frac {\delta}
4\qquad\qquad\qquad (2)$$
\vskip .1in

Apply Erd\"os-Rado's Lemma [C-N] to the family
$\displaystyle{\{F_{\xi}=\cup^n_{i=1}F_{\xi ,i},\xi <\omega_1\}}$ and find 
A uncountable, F finite
such that for $\xi\neq\zeta$ in A
$$F_{\xi}\cap F_{\zeta}=F.$$
\vskip .1in
We set $\lambda_F=\displaystyle{\sum_{\gamma\in F}}\lambda_{\gamma}$ and for $
\xi$ in A
$$\nu_{\xi ,i}=\mu^{\prime}_{\xi ,i}-\frac {d\mu^{\prime}_{\xi ,i}}{
d\lambda_F}$$
\vskip .1in
{\bf Claim:}  For all $i=1,...,n$ the set $B_i=\{\xi\in A:d(T_{\nu_{
\xi ,i}},E_u)\leq\frac {\delta}{16}\}$ is at
most countable.
\vskip .1in
\demo{Proof of the Claim}Suppose that for some $i$ the set $B_i$ is
uncountable.  Then, since $L^1(\lambda_F)$ is separable, there are $
\xi\neq\zeta$ in $B_i$
such that
$$||\frac {d\mu^{\prime}_{\xi ,i}}{d\lambda_F}-\frac {d\mu^{\prime}_{
j,i}}{d\lambda_F}||<\frac {\delta}{16}$$
But then
$$d(T\mu^{\prime}_{\xi ,i}-T\mu^{\prime}_{\zeta ,i},E_u)<\frac {\delta}
4$$
which contradicts inequality (2) and this completes the proof of the
claim.

Choose $\xi_1<\xi$ $_2$ $<...<\xi_n$ in A such that
$$d(T_{\nu_{\xi_i,i}},E_u)>\frac {\delta}{16}\qquad\qquad\qquad\qquad
\qquad\qquad\quad\qquad (3)$$
\vskip .1in
In the rest of the proof we will denote $(\xi_i,i)$ by $\xi_i$.

Notice that the measures $\nu_{\xi_1},...\nu_{\xi_n},\lambda_F$ are pairwise singular.  Choose
$W_1,$..., $W_n$ pairwise disjoint clopen subsets of $\{0,1\}^{\Bbb N}$ pairwise disjoint 
such that for $i=1,...,m$
$$||\nu_{\xi_i}|W^c_i||<\frac {\delta}{128}\text{   and   }||\frac {
d\mu_{\xi_i}}{d_{\lambda_F}}|{n\atop {{{\cup}\atop {j=1}}}}W_j||<\frac {
\delta}{128}\qquad\quad\qquad (4)$$
\vskip .1in
We are ready to prove the desired property.  Indeed, for
$\displaystyle{\lambda_{_i}\geq 0,\text{  }\sum^n_{i=1}\lambda_i=1}$  we have
$$d(\sum^n_{i=1}\lambda_iT\mu^{\prime}_{\xi_i},E_u)\geq d(\sum^n_{
i=1}\lambda_iT\mu^{\prime}_{\xi_i}\restriction{{n\atop {{{\cup}\atop {
j=1}}}}}W_j,E_u)\geq$$
$$d(\sum^n_{i=1}\lambda_i(T\nu_{_{\xi_i}}\restriction W_i),E_u)-\sum^
n_{i
=1}\lambda_i||T\nu_{_{\xi_i}}\restriction{{{\cup}\atop {j\neq i}}}
W_j||-$$
\hskip0.20in\hskip0.20in\hskip0.20in\hskip0.20in\hskip0.20in
\hskip0.20in$\displaystyle{
\sum^n_{i=1}\lambda_i||\frac {d\mu_{\xi_i}}{d_{\lambda_F}}\restriction{{
n\atop {{{\cup}\atop {j=1}}}}}W_j||}$
\vskip .1in
From Lemma 1.2 we get
$$d(\sum^n_{i=1}\lambda_i(T\nu_{_{\xi_i}}\restriction W_i),Eu)\geq\frac 
12\sum^n_{i=1}\lambda_id(T\nu_{_{\xi_i}}\restriction W_i,E_u)$$
and from (3) and (4) we get
$$d(\sum^n_{i=1}\lambda_iT\mu^{\prime}_{\xi_i},E_u)>\frac 12\frac {
3\delta}4-\frac {\delta}{64}=\frac {\delta}{128}$$
Finally from (1) we have
$$d(\sum^n_{i=1}\lambda_ix^{**}_{\xi_i},E_u)>\frac {\delta}{256}$$
So L fails P.C.P., and P.4 ensures that $J^{-1}(L)$ also fails this property.
\qed
\vskip .2in
{\smc 1.4 Remark}   The space E does not contain a subspace isomorphic
to $c_o(\Bbb N)$.  This is because $c_o(\Bbb N)$ 
contains a non RNP closed convex
subset on which norm and weak topologies coincide.  Therefore E
does not embed into a space with an unconditional skipped block
finite dimensional decomposition.  The last follows from the fact that
E fails P.C.P. and it does not contain $c_o(\Bbb N)$.  Finally E semiembeds into
$E_u$ a space with an unconditional basis.
\vskip .1in
\proclaim{1.5 Proposition}The properties RNP and KMP are equivalent
on the subsets of E.  Furthermore if K is closed convex non RNP
subset of E then it contains a subset L 
with a $P\alpha\ell$-representation.\endproclaim
\vskip .1in
\demo{Proof}As we mentioned before if K is closed convex bounded
non RNP then $J$[K] carriers the same properties and it is
contained into $E_u$ which has an unconditional basis.  Therefore, there
exists an L closed convex subset of J[K] with a $P\alpha\ell$-representation
[A-D].   Then $J^{-1}[L]$ has the same property.\qed\enddemo
\vskip .1in
We conclude with the following result.
\vskip .1in
\proclaim{1.6 Theorem}Suppose that X is a separable Banach space
such that $X^{**}/X$ is isomorphic to $\ell^1(\Gamma )$.  
Then X has RNP.\endproclaim
\vskip .1in
\demo{Proof}Assume that X contains a $\delta$-non dentable subset K.  Then
the techniques developed in the proof of Proposition 1.3 shows
that K is non strongly regular.  Actually every
$\displaystyle{\sum^u_{i=1}\lambda_iS_i}$ convex combination of slices will have
diameter greater than $\displaystyle{\frac {\delta}{256}}$.  Hence by a 
result due to 
Bourgain [B], $\ell^1$ embeds into $X^{*}$, and 
by Pelczynski's Theorem [P] M[0,1] 
embeds into $X^{**}$.  But then there exists a sequence $(x^{**}_n
)_{n\in \Bbb N}$ weakly 
convergent to zero and $d(x^{**}_n,X)>\delta$.  This contradicts the Schur 
property of $\ell^1(\Gamma )$.\qed\enddemo
\vskip .1in
{\smc 1.7 Remark}   Odell in [O] has constructed a separable B-space X
with $X^{**}/X\cong\ell^1(2^{\omega}).$  From 
a theorem by Lindenstrauss [L] follows that
every separable B space X and its dual $X^{*}$ are of the form $Z^{
**}/Z$
for some separable Banach space Z.

\settabs\+\hskip1.00in&\hskip1.00in&\hskip1.00in&\hskip1.00in&\hskip1.00in&\cr
\+&&&Department of Mathematics\cr
\+&&&University of Crete\cr
\+&&&Herakleion Crete\cr
\bye